\documentclass[11pt,a4paper]{amsart}
\usepackage{amssymb,amsmath}
\usepackage{mathpazo}

\headsep=1cm \numberwithin{equation}{section}
\newtheorem{theorem}{Theorem}[section]
\newtheorem{lemma}[theorem]{Lemma}

\newtheorem{corollary}[theorem]{Corollary}

\theoremstyle{definition}
\newtheorem{definition}[theorem]{Definition}
\theoremstyle{remark}
\newtheorem{remark}[theorem]{Remark}
\newtheorem{example}[theorem]{Example}

\newcommand{\Ass}{\operatorname{Ass}}
\newcommand{\im}{\operatorname{Im}}

\newcommand{\cd}{\operatorname{cd}}

\newcommand{\V}{\operatorname{V}}

\newcommand{\Ext}{\operatorname{Ext}}
\newcommand{\Supp}{\operatorname{Supp}}
\newcommand{\Tor}{\operatorname{Tor}}
\newcommand{\Hom}{\operatorname{Hom}}

\newcommand{\lo}{\longrightarrow}
\newcommand{\fm}{\frak{m}}
\newcommand{\fp}{\frak{p}}

\newcommand{\fa}{\frak{a}}

\newcommand{\fJ}{\frak{J}}

\begin{document}
\author[Asgharzadeh  and Tousi ]{Mohsen Asgharzadeh and Massoud Tousi}

\title[A unified approach to  local cohomology modules  using Serre classes ]
{A unified approach to  local cohomology modules  \\ using Serre
classes}

\address{M. Asgharzadeh,  Department of Mathematics, Shahid Beheshti University, Tehran, Iran. and Institute for
Studies in Theoretical Physics and Mathematics, P.O. Box 19395-5746,
Tehran, Iran.} \email{asgharzadeh@ipm.ir}
\address{M. Tousi, Department of Mathematics, Shahid Beheshti University,  Tehran, Iran. and Institute for
Studies in Theoretical Physics and Mathematics, P.O. Box 19395-5746,
Tehran, Iran.} \email{mtousi@ipm.ir} 


\begin{abstract} This paper
discusses the connection between the local cohomology modules and
the  Serre classes of $R$-modules. Such connection provided a common
language for expressing some results about the local cohomology
$R$-modules, that has appeared in different papers.
\end{abstract}

\maketitle

\section{Introduction}
Throughout this paper $R$ is a commutative Noetherian  ring, $\fa$
an ideal of $R$ and $M$ an $R$-module.

The proofs of some results concerning local cohomology modules
indicate that these proofs apply to certain subcategories of
$R$-modules that are closed under taking extensions, submodules and
quotients. It should be noted that these kind of subcategories of
$R$-modules are called "Serre classes".  In this paper,
"$\mathcal{S}$" stands for a "Serre class". The aim of the present
paper is to show some results of local cohomology modules are remind
true for any Serre classes. As a general reference for local
cohomology, we refer the reader to the text book \cite{BS}.

Our paper is divided into three sections. In Section 2, we prove the following theorem:\\

\textbf{Theorem A.} Let $ s \in \mathbb{N}_{0}$ and $M$ be an
$R$-module such that $\Ext^{s}_{R}(R/ \fa,M)\in \mathcal{S}$. If
$\Ext^{j}_{R}(R/ \fa,H^{i}_{\fa} (M))\in \mathcal{S}$ for all $i<s$
and  all $j\geq 0$, then $\Hom _{R}(R/ \fa,H^{s}_{\fa} (M))\in
\mathcal{S}$.

One can see that the subcategories of, finitely generated
$R$-modules, minimax $R$-modules, minimax and $\fa$-cofinite
$R$-modules, weakly Laskerian $R$-modules and Matlis reflexive
$R$-modules are examples of Serre classes. So, we can deduce from
Theorem A the main results of \cite{KS}, \cite{BL} and
\cite[Corollary 2.7]{DM}, \cite[Corollary 2.3]{LSY}, \cite[Lemma
2.2]{BN}, \cite[Theorem 1.2]{AKS}, see Corollary 2.4, Corollary 2.5,
Corollary 2.6, Corollary 2.7, Corollary 2.8 and Corollary 2.10.

In Section 3, we investigate the  notation
$\cd_{\mathcal{S}}(\fa,M)$ as the supremum of the integers $i$ such
that $H^{i} _{\fa}(M)\notin \mathcal{S}$. We prove that:

\textbf{Theorem B.} Let  $M$ and $N$ be finitely generated
$R$-modules. Then the following hold:
\begin{enumerate}
\item[i)] Let $t>0$ be an integer. If $N$ has finite Krull dimension and $H_{\fa}^{j}(N)\in\mathcal{S}$ for all $j>t$, then
$H_{\fa}^{t}(N)/\fa H_{\fa}^{t}(N)\in\mathcal{S}$.
\item[ii)] If $\Supp N\subseteq \Supp M$, then $\cd_{\mathcal{S}}(\fa,N)\leq
\cd_{\mathcal{S}}(\fa,M)$.
\end{enumerate}

If $\mathcal{S}$ is equal to the zero class or the class of Artinian
$R$-modules, then we can obtain the results of \cite[Theorem
2.2]{DNT}, \cite[Theorem 2.3]{DY} and \cite[Theorem 3.3]{ADT}.\\

As an application we show that:

\textbf{Theorem C.} Let $M$ be a finitely generated $R$-module. Then
the following hold:

\begin{enumerate}
\item[i)] If $1<d:=\dim M<\infty$, then
$\frac{H^{d-1}_{\fa}(M)}{\fa^n H^{d-1}_{\fa}(M)}$ has finite length
for any $n\in\mathbb{N}$.
\item[ii)] If $(R,\fm)$ is a local ring of Krull
dimension less than 3, then $\Hom_R(R/\fm, H^i _{\fa} (M))$ is a
finitely generated $R$-module for all $i$.
\end{enumerate}

\section{Serre classes and common results on local cohomology modules}
We need the following observation in the sequel.
\begin{lemma}
Let $M\in \mathcal{S}$ and $N$ a finitely generated $R$-module. Then
$\Ext^{j}_{R}(N,M)\in\mathcal{S}$ and
$\Tor^{R}_{j}(N,M)\in\mathcal{S}$ for all $j\geq 0$.
\end{lemma}
{\bf Proof.} We only prove the assertion for the Ext modules and the
proof for the Tor modules is similar. Let $F_{\bullet} : \cdots
\longrightarrow F_{1} \longrightarrow F_0 \longrightarrow 0$ be a
finite free resolution of $N$. If $F_i = R^{n_i}$ for some integer
$n_i$, then $\Ext^i _R(N,M) = H^i(\Hom_R(F_{\bullet},M))$ is a
subquotient of $M^{n_i}$. Since $\mathcal{S}$ is a Serre class, it
follows  that $\Ext^i _R(N,M)\in\mathcal{S}$ for all $i\geq 0$.
$\Box$\\

The following is one of the main result of this section.

\begin{theorem}
Let $ s \in \mathbb{N}_{0}$ and $M$ be an $R$-module such that
$\Ext^{s}_{R}(R/ \fa,M)\in \mathcal{S}$. If $\Ext^{j}_{R}(R/
\fa,H^{i}_{\fa} (M))\in \mathcal{S}$ for all $i<s$ and  all $j\geq
0$, then $\Hom _{R}(R/ \fa,H^{s}_{\fa} (M))\in \mathcal{S}$.
\end{theorem}

{\bf Proof.} We use induction on $s$. From the isomorphism
$\Hom_{R}(\frac{R}{\fa},M) \cong \Hom_{R}(\frac{R}{\fa},
\Gamma_{\fa}(M))$, the case $s=0$ follows. Now suppose inductively
that $s>0$ and that the assertion holds for $s-1$. Let $L=M/
\Gamma_{\fa}(M)$. Then there exists the exact sequence
$$0 \longrightarrow \Gamma_{\fa}(M) \longrightarrow M\longrightarrow
L \longrightarrow 0 .$$ This sequence induces the exact sequences
$$ \Ext^{j}_{R}(R/ \fa,M) \longrightarrow \Ext^{j}_{R}(R/
\fa,L)\longrightarrow \Ext^{j+1}_{R}(R/ \fa,\Gamma_{\fa}(M))$$ for
all $j\geq0$. On the other hand, we have $H^{i}_{\fa} (M)\cong
H^{i}_{\fa} (L)$ for all $i\geq 1$ and $\Gamma_{\fa}(L)=0$.  Also,
by our assumption we have $\Ext^{j}_{R}(R/ \fa,\Gamma_{\fa}(M))\in
\mathcal{S}$ for all $j\geq 0$. Hence we can replace $M$ by $M/
\Gamma_{\fa}(M)$. Therefore $\Gamma_{\fa}(M)=0$. Let $E_{R}(M)$ be
an injective envelope of $M$. Then we have  the exact sequence $$0
\longrightarrow M \longrightarrow E_{R}(M)\longrightarrow N
\longrightarrow 0 .$$ Since
$\Gamma_{\fa}(E_{R}(M))=E_{R}(\Gamma_{\fa}(M))=0$, we have
$H^{i}_{\fa} (N)=H^{i+1}_{\fa} (M)$ for all $i\geq 0$. The fact
$\Hom_{R}(R/ \fa,E_{R}(M))=0$ implies that $\Ext^{j}_{R}(R/
\fa,N)\cong \Ext^{j+1}_{R}(R/ \fa,M)$ for  all $j\geq 0$. So $N$
satisfies our induction hypothesis. Therefore $\Hom _{R}(R/
\fa,H^{s-1}_{\fa} (N))\in \mathcal{S}$. The assertion follows from
$H^{s}_{\fa} (M)\cong H^{s-1}_{\fa} (N)$. $\Box$

\begin{corollary}
Let  assumptions be  as in Theorem 2.2. Let  $N\subseteq H^{s}_{\fa}
(M)$ be such that $\Ext^{1}(R/\fa,N)\in \mathcal{S}$. Then $\Hom
_{R}(R/ \fa,H^{s}_{\fa} (M)/N)\in \mathcal{S}$.
\end{corollary}

{\bf Proof.} The assertion follows from the long $\Ext$ exact
sequence, induced by
$$0\rightarrow N \rightarrow H^{s}_{\fa} (M)\rightarrow H^{s}_{\fa}
(M)/N\rightarrow 0. \  \ \Box$$

The categories of, finitely generated $R$-modules, minimax
$R$-modules \cite[Lemma 2.1]{BN}, weakly Laskerian $R$-modules
\cite[Lemma 2.3]{DM} and Matlis reflexive $R$-modules,  are examples
of Serre classes. Hartshorne defined a module $M$ to be
$\fa$-cofinite, if $\Supp_{R}M \subseteq \V(\fa)$ and $\Ext^{i}
_{R}(R/\fa,M)$ are finitely generated module for all $i$, see
\cite{Har2}. By \cite[Corollary 4.4]{M} the class of $\fa$-cofinite
minimax modules is a Serre class of the category of $R$-modules.
Consequently, we can deduce the following results from Theorem 2.2
and Corollary 2.3.  We denote the set of associated primes of $M$ by
$\Ass_R (M)$. Note that $\Ass_R(\Hom_R(R/\fa,M))=\Ass_R(M)$, for all
$\fa$-torsion
$R$-modules $M$.\\

Khashyarmanesh and Salarian in \cite{KS} proved the following
theorem by the concept of $\fa$-filter regular sequences:

\begin{corollary} Let $M$ be a finitely generated $R$-module and
$t$ an integer. Suppose that  the local cohomology modules $H^{i}
_{\fa}(M)$ are finitely generated for all $i < t$. Then
$\Ass_R(H^{t} _{\fa}(M))$ is finite.
\end{corollary}

On the other hand, Brodmann and Lashgari \cite{BL} generalized this
by the basic homological algebraic methods.

\begin{corollary}  Let $M$ be a finitely generated $R$-module and
$t$ an integer. Suppose that the local cohomology modules $H^{i}
_{\fa}(M)$ are finitely generated for all $i < t$. Then
$\Ass_R(H^{t} _{\fa}(M)/N)$ is finite, for any  finitely generated
submodule $N$ of $H^{i} _{\fa}(M)$.
\end{corollary}

Recall that, from \cite{DM}, an $R$-module $M$ is  weakly Laskerian
R-module, if any quotient of M has a finitely many  associated
primes. Divaani-Aazar and Mafi in \cite[Corollary 2.7]{DM} proved by
the spectral sequences technics:

\begin{corollary} Let $M$ be a weakly Laskerian $R$-module and
$t$ an integer such that $H^{i} _{\fa}(M)$ is weakly Laskerian
modules for all $i < t$. Then $\Ass_R(H^{t} _{\fa}(M))$ is finite.
\end{corollary}

Recall that an $R$-module $M$ is minimax if there is a finitely
generated  submodule $N$ of $M$ such that $M/N$ is Artinian, see
\cite{Z} and \cite{R}.

\begin{corollary} (see \cite[Corollary 2.3]{LSY}). Let $M$ be a minimax
$R$-module. Let $t$ be a non-negative integer such that
$H^{i}_{\fa}(M)$ is a minimax $R$-module for all $i<t$. Let $N$ be a
submodule of $H^{t}_{\fa}(M)$ such that $\Ext^{1}_{R}(R/\fa,N)$ is
minimax. Then $\Hom_R(R/\fa,H^{t}_{\fa}(M)/N))$ is minimax. In
particular $H^{t}_{\fa}(M)/N$ has finitely many associated primes.
\end{corollary}

The following is a key lemma of \cite[Lemma 2.2]{BN}. In fact it is
true without $\fa$-cofinite  condition, see \cite[Theorem 3.2]{BN}.

\begin{corollary} (see \cite[Lemma 2.2]{BN}). Let $M$ be a finitely generated
$R$-module. Let $t$ be a non-negative integer such that
$H^{i}_{\fa}(M)$ are minimax and $\fa$-cofinite $R$-modules for all
$i<t$.  Then $\Hom_{R}(R/\fa,H^{t}_{\fa}(M))$ is  finitely generated
and as a consequence  it has  finitely many associated primes.
\end{corollary}

{\bf Proof.} Set $\mathcal{S}$ be the class of $\fa$-cofinite and
minimax modules. From Theorem 2.2, $\Hom _{R}(R/ \fa,H^{s}_{\fa}
(M))$ is a minimax and $\fa$-cofinite $R$-module. Therefore  we get
that $\Hom_{R}(R/\fa,\Hom _{R}(R/ \fa,H^{s}_{\fa} (M)))\cong \Hom
_{R}(R/ \fa,H^{s}_{\fa} (M))$ is finitely generated. $\Box$

\begin{corollary} Let  $M$ a finitely generated $R$-module
and $\mathcal{S}$ a  Serre class that contains all finitely
generated $R$-modules. Let $t$ be a non-negative integer such that
$H^{i}_{\fa}(M)\in\mathcal{S}$ for all $i<t$.   Then
$\Hom_{R}(R/\fa,H^{t}_{\fa}(M))\in\mathcal{S}$.
\end{corollary}
An immediate consequence of Corollary 2.9 is the following:

\begin{corollary} (see \cite[Theorem 1.2]{AKS}) Let $M$ be a finitely generated
$R$-module. Let $t$ be a non-negative integer such that
$H^{i}_{\fa}(M)$ is finitely generated for all $i<t$.   Then
$\Hom_{R}(R/\fa,H^{t}_{\fa}(M))$ is a finitely generated $R$-module
and as a consequence it has finitely many associated primes.
\end{corollary}

In the proof of Theorem 2.12, we will use the following lemma.

\begin{lemma} Let $(R,\fm)$ be a local ring and $\mathcal{S}$ a
non-zero Serre class. Let $\mathcal{FL}$ be the class of finite
length $R$-modules. Then $\mathcal{FL}\subseteq \mathcal{S}$.
\end{lemma}

{\bf Proof.} Since $\mathcal{S}$ is non-zero, there exists a
non-zero $R$-module $L\in\mathcal{S}$. Let $0\neq m\in L$. Then $Rm
\in\mathcal{S}$. From the natural epimorphism $ Rm\cong
R/(0:_{R}m)\twoheadrightarrow R/\fm$ we obtained that
$R/\fm\in\mathcal{S}$. Let $M\in \mathcal{FL}$ and set
$\ell:=\ell_R(M)$. By induction  on $\ell$, we show  that
$M\in\mathcal{S}$. For the cases $\ell=0,1$, we have nothing to
prove. Now suppose inductively, $\ell>0$ and the result has been
proved for each finite length $R$-module $N$, with $\ell_R(N)\leq
\ell-1$. By definition there is following chain of $R$-submodules of
$M$:
$$0=M_{0}\subseteq M_{1}\subseteq\dots\subseteq M_{\ell}=M$$  such that
$M_{j}/M_{j-1}\cong R/ \fm$. Now the exact sequence
$$0\longrightarrow M_{\ell -1}\longrightarrow M \longrightarrow
R/\fm\longrightarrow0,$$completes the proof. $\Box$

Now we are ready to prove the second main result of this section.

\begin{theorem} Let $(R,\fm)$ be a local ring, $\mathcal{S}$ a
 non-zero Serre class and $M$  a finitely generated $R$-module. Let $t$
be a non-negative integer such that $H^{i}_{\fa}(M)\in\mathcal{S}$
for all $i<t$. Then $\Hom_{R}(R/\fm,H^{t}_{\fa}(M))\in\mathcal{S}$.
\end{theorem}

{\bf Proof.} We do induction on $t$. If $t=0$, then
$\Hom_{R}(R/\fm,H^{0}_{\fa}(M))$ has finite length. So by Lemma
2.11, $\Hom_{R}(R/\fm,H^{0}_{\fa}(M))\in\mathcal{S}$. Now suppose
inductively, $t>0$ and the result has been proved for all integer
smaller than $t$. We have $H^{i}_{\fa}(M)\cong H^{i}_{\fa}(M/\Gamma
_{\fa}(M))$  for all $i>0$. Hence we may assume that $M$ is
$\fa$-torsion free. Take $x\in\fa\setminus \bigcup _{\fp \in
\Ass_{R}M}\fp$.  From the exact sequence
$$0\longrightarrow M \stackrel{x}\longrightarrow M\longrightarrow
M/xM\longrightarrow 0$$ we  deduced the long exact sequence of local
cohomology modules, which shows that
$H^{j}_{\fa}(M/xM)\in\mathcal{S}$ for all $j<t-1$. Thus,
$\Hom_R(R/\fm,H^{t-1}_{\fa}(M/xM))\in\mathcal{S}$.

Now, consider the  long exact sequence $$\cdots\longrightarrow
H^{t-1}_{\fa}(M/xM)\longrightarrow
H^{t}_{\fa}(M)\stackrel{x}\longrightarrow
H^{t}_{\fa}(M)\longrightarrow \cdots,$$ which induces the following
exact sequence
$$0\longrightarrow
H^{t-1}_{\fa}(M)/xH^{t-1}_{\fa}(M)\longrightarrow
H^{t-1}_{\fa}(M/xM)\longrightarrow
(0:_{H^{t}_{\fa}(M)}x)\longrightarrow 0.$$ From this we get the
following exact sequence
$$ \Hom_R(\frac{R}{\fm}, H^{t-1}_{\fa}(\frac{M}{xM}))
\longrightarrow \Hom_R(\frac{R}{\fm}, (0:_{H^{t}_{\fa}(M)}x))
\longrightarrow\Ext^1_R(\frac{R}{\fm},\frac{H^{t-1}_{\fa}(M)}{xH^{t-1}_{\fa}(M)}).$$
By Lemma 2.1, $\Ext^1_R(\frac{R}{\fm},\frac{H^{t-1}_{\fa}(M)}
 {xH^{t-1}_{\fa}(M)})\in\mathcal{S}$. Therefore, $\Hom_R(R/\fm, (0:_{H^{t}_{\fa}(M)}x))\in\mathcal{S}$.
The following completes the proof:
\[\begin{array}{ll}
\Hom_{R}(\frac{R}{\fm},H^{t}_{\fa}(M))&\cong\Hom_{R}(R/\fm\otimes_R
R/xR,H^{t}_{\fa}(M))\\&\cong\Hom_R(R/\fm, (0:_{H^{t}_{\fa}(M)}x)). \
\ \Box
\\
\end{array}\]

\begin{example} In Theorem 2.12, the assumption $\mathcal{S}\neq\{0\}$ is
necessary. To see this, let $(R,\fm)$ be a local Gorenstein ring of
positive dimension $d$. Then $H^{i}_{\fm}(R)=0$ for $i<d$. But
$\Hom_R(R/\fm,H^{d}_{\fm}(R))\cong\Hom_R(R/\fm,E)\cong R/\fm\neq0$,
where $E$ is an injective envelope of $R/\fm$.
\end{example}

As an immediate result of Theorem 2.12 (or Corollary 2.10) we have
the following corollary.

\begin{corollary}Let $(R,\fm)$ be a local ring  and $M$ a finitely generated
$R$-module. Let $t$ be an non-negative integer such that
$H^{i}_{\fa}(M)$ is a finitely generated $R$-module for all $i<t$.
Then $\Hom_{R}(R/\fm,H^{t}_{\fa}(M))$ is a finitely generated
$R$-module.
\end{corollary}

Let $(R,\fm)$ be a local ring. The third of Huneke's four problems
in local cohomology \cite{Hu} is to determine when $H^i _{\fa} (M)$
is Artinian for a finitely generated $R$-module M. The mentioned
problem may be separated into two subproblems:\begin{enumerate}
\item[i)]When is $\Supp_R (H^i _{\fa} (M)) \subseteq \{\fm\}$?
\item[ii)] When is $\Hom_R(R/\fm, H^i _{\fa} (M))$ finitely generated?
\end{enumerate}
Huneke formalized the following conjecture, see \cite[Conjecture
4.3]{Hu}.

\textbf{Conjecture.} Let $(R,\fm)$ be a regular local ring  and
$\fa$ be an ideal of $R$. For all $i$,  $\Hom_R(R/\fm, H^i _{\fa}
(R))$ is finitely generated.

It is known that if R is an unramified regular local ring, then
$\Hom_R(R/\fm, H^i _{\fa} (R))$  is finitely generated, for all $i$
(see \cite{HS}, \cite{L1}, \cite{L2}). The first example of a local
cohomology module with an infinite dimensional socle was given in
\cite{Har2} by Hartshorne. The Hartshorne's famous example is a
three dimensional local ring.

As the first application, the following provides a positive answer
of the conjecture, for all local rings of Krull dimension less than
3.

\begin{corollary} Let $(R,\fm)$ be a local ring  of dimension less
than 3, and $M$ a finitely generated $R$-module. Then $\Hom_R(R/\fm,
H^i _{\fa} (M))$ is a finitely generated $R$-module for all $i$.
\end{corollary}

{\bf Proof.}  First assume that $\dim R=2$. The cases $i=0$ and
$i>2$ are trivial, since $H^0 _{\fa} (M)$ is finitely generated
$R$-module and $H^i _{\fa} (M)=0$ for all $i>2$. Note that $H^2
_{\fa} (M)$ is an Artinian $R$-module. Therefore,  $\Hom_R(R/\fm,
H^2 _{\fa} (M))$ is a finitely generated $R$-module. In the case
$i=1$ one can get the desired result from Corollary 2.14.

If $\dim R\leq 1$ we can obtain the desired result in similar way.
$\Box$

\begin{remark} Let $n$ be an integer  grater than $2$.
Then  \cite[Theorem 1.1]{MV} and the discussion before than
\cite[Question 2.1]{MV}, provided  an $n$-dimensional regular local
ring $(R,\fm)$ and a finitely generated $R$-module $M$ such that
$\Hom_{R}(R/\fm,H^{t}_{\fa}(M))$ is not  finitely generated
$R$-module, for some $t\in \mathbb{N}$ and some ideal
$\fa\vartriangleleft R$.
\end{remark}

\section{Serre cohomological dimension}
In the proof of the following theorem we use the method of the proof
of \cite[Theorem 3.3]{ADT}.

\begin{theorem} Let $\fa$ be an ideal of  $R$ and $M$ a
weakly Laskerian $R$-module of finite Kryll dimension. Let $t>0$ be
an integer. If $H_{\fa}^{j}(M)\in\mathcal{S}$ for all $j>t$, then
$H_{\fa}^{t}(M)/\fa H_{\fa}^{t}(M)\in\mathcal{S}$.

\end{theorem}

{\bf Proof.}  We use induction on $d:=\dim M$. The case $d=0$, is
easy, because  $H^{t}_{\fa}(M)=0$. Now suppose inductively, $\dim
M=d>0$ and the result has been proved for all $R$-modules of
dimension smaller than $d$. We have $H^{i}_{\fa}(M)\cong
H^{i}_{\fa}(M/\Gamma _{\fa}(M))$  for all $i>0$. Also $M/\Gamma
_{\fa}(M)$ has dimension not exceeding $d$. So we may assume that
$M$ is $\fa$-torsion free.  Let $x\in\fa\setminus \bigcup _{\fp \in
\Ass_{R}M}\fp$. Then $M/xM$ is weakly Laskerian and $\dim M/xM\leq
d-1$. The exact sequence
$$0\longrightarrow M \stackrel{x}\longrightarrow M\longrightarrow
M/xM\longrightarrow 0$$ induces the long exact sequence of local
cohomology modules, which shows that
$H^{j}_{\fa}(M/xM)\in\mathcal{S}$ for all $j>t$. By induction
hypothesis $H^{t}_{\fa}(M/xM)/\fa H^{t}_{\fa}(M/xM)\in\mathcal{S}$.

Now, consider the exact sequence $$
H^{t}_{\fa}(M)\stackrel{x}\longrightarrow
H^{t}_{\fa}(M)\stackrel{f}\longrightarrow
H^{t}_{\fa}(M/xM)\stackrel{g}\longrightarrow H^{t+1}_{\fa}(M),$$
which induces the following two exact sequences
$$
H^{t}_{\fa}(M)\stackrel{x}\longrightarrow
H^{t}_{\fa}(M)\longrightarrow \im f\longrightarrow 0,$$
$$0\longrightarrow \im f \longrightarrow H^{t}_{\fa}(M/xM)
\longrightarrow\im g \longrightarrow 0.$$ Therefore we can obtain
the following two exact sequences:
$$H^{t}_{\fa}(M)/\fa H^{t}_{\fa}(M)\stackrel{x}\longrightarrow
H^{t}_{\fa}(M)/\fa H^{t}_{\fa}(M) \longrightarrow \im f/\fa \im
f\longrightarrow0,$$
$$Tor^{R}_{1}(R/\fa,\im g )\lo \im f/\fa\im f\lo
H^{t}_{\fa}(M/xM)/\fa H^{t}_{\fa}(M/xM)\lo \im g/\fa \im g\lo 0.$$
Since $x\in \fa$, from a preceding exact sequence, we get that $\im
f/\fa\im f\cong H^{t}_{\fa}(M)/\fa H^{t}_{\fa}(M)$. By Lemma 2.1, we
have $Tor^{R}_{1}(R/\fa,\im g)\in\mathcal{S}$. Also
$H^{t}_{\fa}(M/xM)/\fa H^{t}_{\fa}(M/xM)\in\mathcal{S}$. So $\im
f/\fa\im f\in\mathcal{S}$. Now the claim follows. $\Box$\\

The second of our applications is the following corollary.

\begin{corollary} Let
$M$ be a finitely generated $R$-module of finite Krull dimension
$d>1$. Then $\frac{H^{d-1}_{\fa}(M)}{\fa^n H^{d-1}_{\fa}(M)}$ has
finite length for any $n\in\mathbb{N}$.
\end{corollary}

{\bf Proof.} We have $H^{d-1}_{\fa}(M)=H^{d-1}_{\fa^n}(M)$. So it is
enough to prove the desired result for $n=1$. By \cite[Proposition
5.1]{M}, $H^d _\fa (M)$ is $\fa$-cofinite and Artinian. Set
$\mathcal{S}:=\{N:N \textit{ is a $\fa$-cofinite and minimax
$R$-module}\}$. In view of Theorem 3.1, we get that the $R$-module
$\frac{H^{d-1}_{\fa}(M)}{\fa H^{d-1}_{\fa}(M)}$ is $\fa$-cofinite.
So $\Hom_R(R/\fa,\frac{H^{d-1}_{\fa}(M)}{\fa H^{d-1}_{\fa}(M)})\cong
\frac{H^{d-1}_{\fa}(M)}{\fa H^{d-1}_{\fa}(M)}$ is a finitely
generated $R$-module. Set $\mathcal{S}:=\{N:N \textit{ is an
Artinian $R$-module}\}$. Again by Theorem 3.1, we get that the
$R$-module $\frac{H^{d-1}_{\fa}(M)}{\fa H^{d-1}_{\fa}(M)}$ is an
Artinian $R$-module. Consequently, the $R$-module
$\frac{H^{d-1}_{\fa}(M)}{\fa H^{d-1}_{\fa}(M)}$  has finite length.
$\Box$

\begin{example}
In Corollary 3.2,  if $t<\dim M-1$, then it can be seen that
$H_{\fa}^{t}(N)/\fa H_{\fa}^{t}(N)$ is not necessarily has finite
length. To see this, let $R:=k[[X_{1},\cdots,X_{4}]]$,
$\fJ_{1}:=(X_{1},X_{2}), \fJ_{2}:=(X_{3},X_{4})$ and
$\fa:=\fJ_{1}\cap \fJ_{2}$, where $k$ is a field. By Mayer-Vietoris
exact sequence we get that $H_{\fa}^{2}(R)\cong
H_{\fJ_{1}}^{2}(R)\oplus H_{\fJ_{2}}^{2}(R)$. Now consider the
following isomorphisms

\[\begin{array}{ll}
H_{\fa}^{2}(R)/\fa H_{\fa}^{2}(R)&\cong (H_{\fJ_{1}}^{2}(R)/\fa
H_{\fJ_{1}}^{2}(R))\oplus (H_{\fJ_{2}}^{2}(R)/\fa
H_{\fJ_{2}}^{2}(R))\\&\cong H_{\fJ_{1}}^{2}(R/\fa)\oplus
H_{\fJ_{2}}^{2}(R/\fa).
\\
\end{array}\]

By Hartshorne-Lichtenbaum vanishing Theorem,
$H_{\fJ_{1}}^{2}(R/\fa)\neq 0$. Therefore the cohomological
dimension of $R/\fa $ with respect to $\fJ_{1}$ is two. By \cite[
Remark 2.5]{Hel} the local cohomology $H_{\fJ_{1}}^{2}(R/\fa)$ is
not finitely generated. Consequently, $H_{\fa}^{2}(R)/\fa
H_{\fa}^{2}(R)$ is not finitely generated.
\end{example}

\begin{definition}
Let $M$ be an $R$-module and $\fa$ an ideal of $R$. For a Serre
class $\mathcal{S}$, we define $\mathcal{S}$-cohomological dimension
of $M$ with respect to $\fa$, by
$\cd_{\mathcal{S}}(\fa,M):=\sup\{i\in\mathbb{N}_0:H_{\fa}^{i}(M)\notin\mathcal{S}\}$.
\end{definition}

\begin{theorem}
Let  $M$ and $N$ be finitely generated $R$-modules such that
$\Supp_R N\subseteq \Supp_R M$. Then $\cd_{\mathcal{S}}(\fa,N)\leq
\cd_{\mathcal{S}}(\fa,M)$.
\end{theorem}

{\bf Proof.} It is enough to show that if
$i>\cd_{\mathcal{S}}(\fa,M)$, then $H^{i} _{\fa}(N)\in\mathcal{S}$.
We prove this by descending induction on $i$ with
$\cd_{\mathcal{S}}(\fa,M)< i \leq dim (M)+1$. Note that any non
empty Serre class containing the zero module. By Grothendieck's
vanishing theorem, in the case $i = \dim M +1$ we have nothing to
prove. Now suppose $\cd_{\mathcal{S}}(\fa,M)< i \leq \dim M$ and we
have proved that $H_{\fa}^{i+1}(K)\in\mathcal{S}$ for each finitely
generated $R$-module $K$ with $\Supp_{R}K\subseteq\Supp_{R}M$. By
theorem of Gruson \cite[Theorem 4.1]{V}, there is a chain
$$0=N_{0}\subset N_{1}\subset\dots\subset N_{\ell}=N$$ such that each
of the factors $N_{j}/N_{j-1}$ is a homomorphic image of a direct
sum of finitely many copies of $M$. By using short exact sequences,
the situation can be reduced to the case $\ell = 1$. Therefore, for
some positive integer $n$ and some finitely generated $R$-module
$L$, there exists an exact sequence $0\longrightarrow
L\longrightarrow M^{n}\longrightarrow N\longrightarrow 0$. Thus we
have the following long exact sequence
$$\cdots\longrightarrow H^{i}
_{\fa}(L) \longrightarrow H^{i}_{a}(M^{n})\longrightarrow
H^{i}_{\fa}(N)\longrightarrow H^{i+1} _{\fa}(L)\longrightarrow\cdots
.$$ By the inductive assumption  $H^{i+1}_{\fa}(L)\in \mathcal{S}$.
Since $H^{i} _{\fa}(M^{n})\in\mathcal{S}$ we get that
$H^{i}_{\fa}(N)\in\mathcal{S}$. This completes the inductive step.
$\Box$

Let $\mathcal{A}$ be the class of Artinin $R$-modules. Recall that
in the literatures the notion $\cd_{\{0\}}(\fa,M)$ is denote by
$\cd(\fa,M)$ and $\cd_{\mathcal{A}}(\fa,M)$ by $q _{\fa}(M)$. Here,
we record several immediate consequences of Theorem 3.5.

\begin{corollary}(see \cite[Theorem 2.2]{DNT})
Let $M$ and $N$ be finitely generated $R$-modules such that $\Supp_R
N \subseteq \Supp_R M $. Then $\cd(\fa,N)\leq \cd (\fa,M)$.
\end{corollary}

\begin{corollary}
Let  $M$  be a finitely generated $R$-module. Then
$\cd_{\mathcal{S}}(\fa,M)=\max\{\cd_{\mathcal{S}}(\fa,R/\fp):\fp\in\Ass_{R}
M\}$.
\end{corollary}

{\bf Proof.} Let $N:=\oplus_{\fp \in \Ass_RM}R/\fp$. Then $N$ is
finitely generated and $\Supp_RN=\Supp_RM$. In view of  Theorem 3.5,
$\cd_{\mathcal{S}}(\fa,M)=\cd_{\mathcal{S}}(\fa,N)=\max\{\cd_{\mathcal{S}}(\fa,R/\fp):\fp\in\Ass_{R}
M\}$. $\Box$

\begin{corollary}(see \cite[Theorem 2.3]{DY})
Let  $M$ and $N$ be finitely generated $R$-modules such that
$\Supp_R N \subseteq \Supp_RM$. Then $q_{\fa}(N)\leq q _{\fa}(M)$.
\end{corollary}

We denote by $q(\fa)$ the supremum of all integers $j$ for which
there is a finitely generated $R$-module $M$, with $H_{\fa}^j(M)$
not Artinian. It was proved by Hartshorn \cite{Har1} that $q(\fa)$
is the supremum of all integers $j$ for which $H_{\fa}^j(R)$ is not
Artinian. The following is a generalization of this result.

\begin{corollary}
$\cd_{\mathcal{S}}(\fa,R) = \sup \{\cd_{\mathcal{S}}(\fa,N)|N
\textit{ is a finitely generated R-module}\}$. In particular  if
$H^{j}_{\fa}(R)\in\mathcal{S}$ for all $j>\ell$, then
$H^{j}_{\fa}(M)\in\mathcal{S}$ for all $j>\ell$ and all finitely
generated $R$-module $M$.
\end{corollary}

%%%%%%%%%%%%%%%%%%%%%%%%%%%%%%%%%%%%%%%%%%%%%%%%%%%


\begin{thebibliography}{99}

\bibitem[AKS]{AKS}{J. Asadollahi}, { K. Khashyarmanesh}, { S. Salaarian}, {\it On the finiteness properties
of the generalized local cohomology modules},  Comm. Alg., {\bf 30},
(2), (2002), 859--867.

\bibitem[ADT]{ADT}{ M. Asgharzadeh}, {K. Divaani-Aazar }, {M. Tousi }, {\it
Finiteness dimension of local cohomology modules and its dual
notion}, arXiv: 0803.0461v1.

\bibitem[BN]{BN}{ K. Bahmanpour}, {R. Naghipour}, {\it
On the cofiniteness of local cohohomology modules}, Proc. Amer.
Math. Soc., to appear.

\bibitem[BL]{BL}{ M. Brodmann}, { A. Lashgari},
{\it A finiteness result for associated primes of local cohomology
modules} Proc. Amer. Math. Soc.,   {\bf 128}, (10), (2000),
2851--2853.

\bibitem[BS]{BS}{ M. Brodmann}, { R. Y. Sharp}, {\it
Local cohomology, An Algebraic Introduction with Geometric
Applications}, Cambridge Univ. Press , {\bf 60}, (1998).

\bibitem[DM]{DM}{ K. Divaani-Aazar}, { A. Mafi}, {\it Associated prime
of local cohomology modules}, Proc. Amer. Math. Soc.,  {\bf133},
(3), (2005), 655--660.

\bibitem[DNT]{DNT}{ K. Divaani-Aazar}, { R. Naghipour}, { M. Tousi}, {\it Cohomological dimension
of certain algebraic varieties}, Proc. Amer. Math. Soc., {\bf130},
(12), (2002), 3537--3544.

\bibitem[DY]{DY}{ M. T. Dibaei}, { S. Yassemi }, {\it Associated primes and cofiniteness of
local cohomology modules}, Manuscripta Math., {\bf117}, (2), (2005),
199--205.

\bibitem[Har1]{Har1} {R. Hartshorne}, {\it Cohomological dimension of algebraic varieties},
Ann. of Math., {\bf 88} (1968), 403--450.

\bibitem[Har2]{Har2}{R. Hartshorne},  {\it Affine Duality and
Cofiniteness},  Invent. Math. 9, (1970), 145--164.

\bibitem[Hel]{Hel}{Hellus}, {\it A Note on the injective dimension of local
cohomology modules}, Proc., AMS, {\bf136},  (2008), 2313--2321.

\bibitem[HS]{HS}{C. Huneke}, {R. Sharp}, {\it Bass numbers of local cohomology
modules}, Trans. American Math. Soc. 339 (1993), 765-779.

\bibitem[Hu]{Hu}{C. Huneke}, {\it Problems on local cohomology, Free Resolutions
in commutative algebra and algebraic geometry (Sundance, Utah,
1990)}, Research Notes in Mathematics 2, Boston, MA, Jones and
Bartlett Publishers, (1994), 93--108.

\bibitem[KS]{KS}{ K. Khashyarmanesh}, { S. Salaarian}, {\it
On the associated primes of local cohomology modules},  Comm. Alg.
{\bf 27}, (1999), 6191-6198.

\bibitem[LSY]{LSY}{ B. Lorestani}, {P. Sahandi}, { S. Yassemi}, {\it
Artinian local cohomology modules},  To appear in Canadian
Mathematical Bulletin.

\bibitem[L1]{L1}{G. Lyubeznik}, {\it Finiteness properties of local cohomology modules
(an application of D- modules to commutative algebra)}, Inv. Math.
113 (1993), 41-55.

\bibitem[L2]{L2}{G. Lyubeznik}, {\it Finiteness properties of local
cohomology modules for regular local rings of mixed characteristic:
the unramified case}, Comm. Alg., {\bf 28} (12),  (2000),
5867--5882.

\bibitem[M]{M}{ L. Melkersson},  {\it Modules cofinite with respect to an
ideal}, Journal of Algebra {\bf 285}, (2005), 649-668.

\bibitem[MV]{MV}{T.  Marley}, {J. C. Vassilev}, {\it
Local cohomology modules with infinite dimensional socles},  Proc.,
AMS, {\bf132}, (12), (2004), 3485--3490.

\bibitem[R]{R}{ P. Rudlof},  {\it
On minimax and related modules}, Canada J. math. {\bf 44}, (1992),
154-166.

\bibitem[V]{V} {W. V. Vasconcelos}, {\it Divisor theory in module categories},
North-Holland Mathematics Studies, {\bf 14}.

\bibitem[Z]{Z} {H. Zochinger}, {\it Minimax  modules},
Journal of Algebra {\bf 102}, (1986), 1-32.

\end{thebibliography}
\end{document}